\documentclass[10pt]{article}

\usepackage{latexsym}
\usepackage{amsfonts}
\usepackage{amsmath}
\usepackage{mathrsfs}  
\usepackage{dsfont}    
\usepackage{bbold}     
\usepackage[english]{babel}
\usepackage{caption}
\usepackage{epsfig}
\usepackage{float}
\usepackage{isolatin1}
\usepackage{subfigure}

\DeclareMathOperator{\Val}{Val}

\newtheorem{theorem}{Theorem}

\newtheorem{lemma}[theorem]{Lemma}

\newcommand{\zerarcounters}{\setcounter{equation}{0}\setcounter{theorem}{0}}

\newcommand{\ZZZ}{\mathds{Z}}
\newcommand{\CCC}{\mathds{C}}
\newcommand{\NNN}{\mathds{N}}

\newcommand{\RRR}{\mathds{R}}

\newcommand{\MM}{{\mathcal M}}

\newcommand{\RR}{{\mathcal R}}
\newcommand{\SSS}{{\mathcal S}}
\newcommand{\TT}{{\mathcal T}}

\newcommand{\VV}{{\mathcal V}}

\newcommand{\RRRR}{{\mathfrak R}}

\newcommand{\Fullbox}{{\rule{2.0mm}{2.0mm}}}

\newcommand{\EP}{\hfill\Fullbox\vspace{0.2cm}}
\newcommand{\prova}{\noindent{\it Proof. }}
\newcommand{\io}{\infty}
\newcommand{\e}{\varepsilon}
\newcommand{\al}{\alpha}

\newcommand{\n}{\nu}

\newcommand{\x}{\xi}

\newcommand{\ka}{\kappa}

\newcommand{\om}{\omega}

\newcommand{\laa}{\langle}
\newcommand{\raa}{\rangle}

\newcommand{\oo}{\boldsymbol{\omega}}
\newcommand{\nn}{\boldsymbol{\nu}}
\newcommand{\pps}{\boldsymbol{\psi}}
\newcommand{\vzero}{\boldsymbol{0}}

\def\ins#1#2#3{\vbox to0pt{\kern-#2 \hbox{\kern#1 #3}\vss}\nointerlineskip}

\begin{document}

\title{\bf Quasi-periodic attractors, Borel summability and\\
the Bryuno condition for strongly dissipative systems}

\author
{\bf Guido Gentile$^\dagger$, Michele V. Bartuccelli$^\ast$ and 
Jonathan H.B. Deane$^\ast$ 
\vspace{2mm} \\ \small 
$^\dagger$Dipartimento di Matematica, Università di Roma Tre, Roma,
I-00146, Italy.
\\ \small 
E-mail: gentile@mat.uniroma3.it
\\ \small
$^\ast$Department of Mathematics and Statistics,
University of Surrey, Guildford, GU2 7XH, UK.
\\ \small 
E-mails:
m.bartuccelli@surrey.ac.uk, j.deane@surrey.ac.uk}

\date{}

\maketitle

\begin{abstract}
We consider a class of ordinary differential equations describing
one-dimensional analytic systems with a quasi-periodic forcing term
and in the presence of damping. In the limit of large damping,
under some generic non-degeneracy condition on the force, there are
quasi-periodic solutions which have the same frequency vector as the
forcing term. We prove that such solutions are Borel summable at the
origin when the frequency vector is either any one-dimensional
number or a two-dimensional vector such that the ratio of its
components is an irrational number of constant type.
In the first case the proof given simplifies that provided
in a previous work of ours. We also show that in any dimension $d$,
for the existence of a quasi-periodic solution with the same
frequency vector as the forcing term, the standard Diophantine
condition can be weakened into the Bryuno condition.
In all cases, under a suitable positivity condition, the
quasi-periodic solution is proved to describe a local attractor.
\end{abstract}




\zerarcounters
\section{Introduction}
\label{sec:1}

In this paper we pursue the study started in \cite{GBD,BDG}.
We consider one-dimensional systems with a quasi-periodic forcing
term in the presence of strong damping, described
by ordinary differential equations of the form
\begin{equation}
\e \ddot x + \dot x + \e g(x) = \e f(\oo t) ,
\label{eq:1.1} \end{equation}
where $\oo\in\RRR^{d}$ is the frequency vector,
$g(x)$ and $f(\pps)$ are functions analytic in their arguments,
with $f$ quasi-periodic, i.e.
\begin{equation}
f(\pps) = \sum_{\nn\in \ZZZ^{d}}
{\rm e}^{i\nn \cdot \pps} f_{\nn} ,
\label{eq:1.2} \end{equation}
with average $\laa f \raa = f_{\vzero}$, and $\e>0$ is a real parameter,
physically representing the inverse of the damping coefficient.
With $\cdot$ we are denoting the scalar product in $\RRR^{d}$.
A Diophantine condition is assumed on $\oo$ for $d>1$, that is
\begin{equation} \left| \oo\cdot\nn \right| \ge C_{0} |\nn|^{-\tau} \qquad
\forall \nn\in\ZZZ^{d} \setminus\{\vzero\} ,
\label{eq:1.3} \end{equation}
where $|\nn|=|\nn|_{1}\equiv|\n_{1}|+\ldots+|\n_{d}|$,
and $C_{0}$ and $\tau$ are positive constants.
The set of vectors satisfying the condition (\ref{eq:1.3})
is non-void for $\tau\ge d-1$ and is of full measure for $\tau>d-1$.
For $d=1$ we denote the vectors without boldface; in that case
$\om$ will be called the frequency number.

In \cite{GBD} we showed that, under the non-degeneracy condition
\begin{equation} \exists c_{0}\in\RRR \hbox{ such that } g(c_{0}) =
f_{\vzero} \hbox{ and }
g'(c_{0}) \neq 0 ,
\label{eq:1.4} \end{equation}
the system (\ref{eq:1.1}) admits a quasi-periodic solution
$x(t;\e)$ with the same frequency vector as the forcing.
Such a solution can be obtained by a suitable summation
of the formal power series
\begin{equation} x_{0}(t;\e) := \sum_{k=0}^{\io} \e^{k} x^{(k)}(t) ,
\qquad x^{(k)}(t) = \sum_{\nn\in\ZZZ^{d}}
{\rm e}^{i\nn\cdot\oo t} x^{(k)}_{\nn} ,
\label{eq:1.5} \end{equation}
which solves the equations of motion order by order. For $d=1$
(periodic forcing) the series (\ref{eq:1.4}) is Borel summable in $\e$.
In \cite{BDG} we also showed that if $g'(c_{0}) >0$, for any $d$
such a solution is locally an attractor. In some cases, for instance if
$g(x)=x^{2p+1}$, $p\in\NNN$, and $f_{\vzero}>0$, the attractor is global.

In this paper we first give a different (simpler) proof of
Borel summability in the periodic case (Section \ref{sec:2}),
then we prove that the formal series for the solution turns out to be
Borel summable also for $d=2$ and $\tau=1$ (Section \ref{sec:3});
this corresponds to frequency vectors with components such that their
ratios are irrational numbers of constant type (i.e. numbers with
bounded partial quotients in their continued fraction expansion).
The proof does not rely on Nevanlinna-type theorems \cite{S},
but consists in checking directly that the
conditions for the formal series of the solution to be Borel summable
are satisfied, and follows the same strategy introduced in \cite{CGGG}
to investigate Borel summability of lower-dimensional tori.

Finally in Section \ref{sec:4} we show how to relax the
Diophantine condition. We show that, in order to have the
same results on existence and attractivity of the
quasi-periodic solution, one can take $\oo$ to be a
Bryuno vector, that is one can assume that, by defining
\begin{equation}
B(\oo) = \sum_{n=0}^{\infty} \frac{1}{2^{n}} \log
\frac{1}{\al_{n}(\oo)} , \qquad \al_{n}(\oo) = \inf_{|\nn|\le 2^{n}}
|\oo\cdot\nn| ,
\label{eq:1.6} \end{equation}
then $\oo$ satisfies the Bryuno condition $B(\oo)<\infty$.
More formal statements will be given in next sections.

\zerarcounters
\section{Borel summability for $\boldsymbol{d}
\boldsymbol{=}\boldsymbol{1}$}
\label{sec:2}

First of all let us recall the definition of Borel summability \cite{S}.
Let $f(\e)=\sum_{n=1}^{\infty} a_{n}\e^{n}$ a formal power series (which
means that the sequence $\{a_{n}\}_{n=1}^{\infty}$ is well defined).
We say that $f(\e)$ is \textit{Borel summable} if
\begin{enumerate}
\item $B(p):=\sum_{n=1}^{\infty} a_{n} p^{n}/n!$
converges in some circle $|p|<\delta$,
\item $B(p)$ has an analytic continuation to a neighbourhood
of the positive real axis, and
\item $g(\e)=\int_{0}^{\infty} {\rm e}^{-p/\e} B(p) \, {\rm d}p$
converges for some $\e>0$.
\end{enumerate}
Then the function $B(p)$ is called the \textit{Borel transform} of $f(\e)$,
and $g(\e)$ is the \textit{Borel sum} of $f(\e)$. Moreover if the
integral defining $g(\e)$ converges for some $\e_{0}>0$ then it
converges in the circle $\hbox{Re}\,\e^{-1} > \hbox{Re}\,\e_{0}^{-1}$.
A function which admits the formal power series expansion $f(\e)$
is called Borel summable if $f(\e)$ is Borel summable;
in that case the function equals the Borel sum $g(\e)$.

\begin{theorem}\label{thm:1}
Consider the system (\ref{eq:1.1}) for $d=1$, and assume that
the non-degeneracy condition (\ref{eq:1.4}) is fulfilled.
There exists $\e_{0}>0$ such that for $|\e|<\e_{0}$ there
is a periodic solution $x(t;\e)$ which has the same frequency
number as the forcing term and is Borel summable in $\e$ at the origin.
If $g'(c_{0})>0$ such a solution describes a local attractor.
\end{theorem}

\prova We consider explicitly the case $g(x)=x^{2}$ in (\ref{eq:1.1}),
which corresponds to the varactor equation extensively studied
in \cite{GBD,BDGM,BDG}; the general case can be easily
dealt with by reasoning as in Section VII of \cite{GBD}.
In \cite{GBD} we proved that the formal power series (\ref{eq:1.5})
is well defined and that to any order $k$ one has
\begin{equation}
\left| x^{(k)}_{\n} \right| \le A_{1} \e_{2}^{-k}k! , \qquad
\left| x^{(k)}(t) \right| \le A_{1} \e_{2}^{-k}k! ,
\label{eq:2.1} \end{equation}
for suitable constants $A_{1}$ and $\e_{2}$ (cf. formula (4.5)
in \cite{GBD}). This means that the first condition, in the
definition of Borel summability, is satisfied, with $\delta=\e_{2}$. 

In \cite{GBD} we also proved that the formal power series
can be summed, and gives a function 
\begin{equation}
x(t;\e) = \sum_{k=0}^{\io}
\sum_{\n\in\ZZZ} {\rm e}^{i\om\n t} x^{[k]}_{\n} ,
\label{eq:2.2} \end{equation}
which is real-analytic and periodic in $t$, and analytic in $\e$
in a suitable domain tangent to the imaginary axis at the origin.
The coefficients $x^{[k]}_{\n}$ can be written as
\begin{equation}
x^{[k]}_{\n} = \sum_{\theta \in \TT_{k,\n}}
\Val(\theta) , \qquad
\Val(\theta) = \Big(\prod_{\ell \in L(\theta)} g_{\ell} \Big) 
\Big( \prod_{v \in E(\theta) \cup V(\theta)} F_{v} \Big) ,
\label{eq:2.3} \end{equation}
where the symbols are defined as in Section V of \cite{GBD}.
We briefly recall the basic definitions and notations,
with the purpose of making self-consistent the discussion;
reference should be made to \cite{GBD} for further details.

A tree $\theta$ is a graph, that is a connected set of points and lines,
with no cycle, such that all the lines are oriented toward a unique
point (root) which has only one incident line (root line).
All the points in a tree except the root are denoted nodes.
The orientation of the lines in a tree induces a partial ordering 
relation ($\preceq$) between the nodes. Given two nodes $v$ and $w$,
we shall write $w \preceq v$ every time $v$ is along the path
(of lines) which connects $w$ to the root. We call $E(\theta)$
the set of endpoints in $\theta$, that is the nodes which have no
entering line. The endpoints can be represented either as white bullets
or as black bullets; we denote with $E_W(\theta)$ and $E_{B}(\theta)$
the set of white bullets and the set of black bullets, respectively.
With each endpoint $v$ we associate a mode label $\n_v\in\ZZZ$,
such that $\n_{v}= 0$ if $v \in E_{W}(\theta)$ and $\n_{v}\neq0$
if $v \in E_{B}(\theta)$. We denote with $L(\theta)$ the set of lines
in $\theta$. Since $\ell$ is uniquely identified with the point $v$
which it leaves, we may write $\ell = \ell_v$.  With each line $\ell$
we associate a {\rm momentum} label $\n_{\ell} \in \ZZZ$.
The modes of the endpoints and and the momenta of the lines
are related as follows: if $\ell = \ell_{v}$ one has
\begin{equation}
\n_{\ell} = \sum_{i=1}^{s_{v}} \n_{\ell_{i}} =
\sum_{w \in E_{B}(\theta) : w \preceq v} \n_{w} ,
\label{eq:2.4} \end{equation}
where $s_{v}$ denotes the number of lines entering $v$
(one has $s_{v}=2$ if $g(x)=x^{2}$ in (\ref{eq:1.1}), otherwise
$s_{v} \ge 2$), and $\ell_{1},\ldots,\ell_{s_{v}}$ are the lines
entering $v$. We denote by $V(\theta)$ the set of vertices in $\theta$,
that is the set of points which have at least one entering line.
We set $V_{0}(\theta)=\{v\in V(\theta) : \n_{\ell_{v}}=0\}$.
We call {\it equivalent} two trees which can be transformed into
each other by continuously deforming the lines in such a way that
they do not cross each other. Let $\TT_{k,\n}$ be the set of
inequivalent trees of order $k$ and total momentum $\n$, that is
the set of inequivalent trees $\theta$ such that 
$|V(\theta)|+|E_{B}(\theta)|=k$
and the momentum of the root line is $\n$. We associate
with each line $\ell$ a \textit{propagator}
\begin{equation}
g_{\ell} = \begin{cases}
1/((i\om\n_{\ell})(1+i\e\om\n_{\ell})) , & \n_{\ell} \neq 0, \\
1 , & \n_{\ell}=0 , \end{cases}
\label{eq:2.5} \end{equation}
with each vertex $v$ a node factor
\begin{equation}
F_{v} = \begin{cases} -\e , & v \notin V_{0}(\theta) , \\
-1/2c_{0} , & v \in V_{0}(\theta) , \end{cases}
\label{eq:2.6} \end{equation}
and with each endpoint $v$ a node factor
\begin{equation}
F_{v} = \begin{cases} c_{0} , & v \in E_{W}(\theta) , \\ 
\e f_{\n_{v}} , & v \in E_{B}(\theta) . \end{cases}
\label{eq:2.7} \end{equation}
Then (\ref{eq:2.3}) says that each coefficient $x^{[k]}_{\n}$
is given by the sum over all trees of order $k$ and total
momentum $\n$ of the corresponding values.

It is more convenient to slightly change the definition
of node factors and propagators, by associating the factor $\e$
with the propagator $g_{\ell}$ of the line $\ell$ coming out
from $v$ and not with $v$ itself. In this way the
propagator of any line with $\ell$ momentum $\n_{\ell}\neq0$ is
\begin{equation} g_{\ell} = g(\om\n_{\ell};\e) , \qquad
g(x;\e) = \frac{\e}{ix(1+i\e x)} ,
\label{eq:2.8} \end{equation}
and the only dependence on $\e$ in $\Val(\theta)$ is through
the product of propagators with non-vanishing momentum.\footnote{
Note that $g(x;\e)$ in (\ref{eq:2.8}) has a completely different
meaning with respect to the function $g(x)$ appearing in (\ref{eq:1.1}).
The same \textit{caveat} applies to the propagators $g^{[n]}(x;\e)$
in Section \ref{sec:3}.}

The function (\ref{eq:2.8}) is Borel summable, and its Borel
transform is easily computed to be
\begin{equation}
g_{B}(x;p) = \frac{{\rm e}^{-ipx}}{ix} \qquad \Longrightarrow
\qquad \left| g_{B}(x;p) \right| \le
\frac{{\rm e}^{|{\rm Im} \, p|\,|x|}}{|x|} .
\label{eq:2.9} \end{equation}
Moreover $g_{B}(x;p)$ is an entire function in $p$, and the integral
$ \int_{0}^{\infty} {\rm e}^{-p/\e} g_{B}(x;p) \, {\rm d}p$
converges (absolutely) for all $\e>0$.

For any tree $\theta\in \TT_{k,\n}$ the Borel trasform of $\Val(\theta)$
is given by a constant times the Borel transform of the
product of the propagators with non-zero momentum. One has
\begin{equation}
\left( \Val(\theta) \right)_{B} (p) =
\Big( \prod_{\ell \in L_{0}(\theta)} g_{\ell} \Big) 
\Big( \prod_{v \in E(\theta) \cup V(\theta)} F_{v} \Big)
\Big( \Big( \prod_{\ell \in L_{2}(\theta)}
g_{\ell} \Big)_{B}(p) \Big) ,
\label{eq:2.10} \end{equation}
where we have called $L_{0}(\theta)$ is the set of lines in $L(\theta)$
with zero momentum, and we have set $L_{2}(\theta)=L(\theta)
\setminus L_{0}(\theta)$ (cf. Section IV of \cite{GBD}).
The Borel transform appearing in (\ref{eq:2.10}) equals
the convolution of the Borel transforms of the propagators
with non-zero momentum, so that it can be bounded as
\begin{equation}
\Big| \Big( \prod_{\ell \in L_{2}(\theta)}
g_{\ell} \Big)_{B} (p)\Big| \le
{\mathop{\prod}_{\ell \in L_{2}(\theta)}}^{\!\!\!\!*}
\left| g_{B}(\om\n_{\ell};p)  \right| \le
\Big( \prod_{\ell \in L_{2}(\theta)}
\frac{1}{|\om\n_{\ell}|} \Big) \frac{|p|^{k-1}}{(k-1)!}
\exp \left( |\hbox{Im} \, p| \max_{\ell\in
L_{2}(\theta)} |\om\n_{\ell}| \right) ,
\label{eq:2.11} \end{equation}
where $\prod^{*}$ denotes the convolution product, and
$|\om| < |\om\n_{\ell}|< |\om| \sum_{v\in E_{B}(\theta)}|\n_{v}|$;
cf. Remarks (4) to (6) after Definition 1 in \cite{CGGG} for
properties of the Borel transforms we are using here.

Therefore, for $p$ in any strip $\Sigma_{\sigma}=\{p\in\CCC:
|{\rm Im}\,p| < \sigma\}$ of the real axis, we have
\begin{equation}
\Big| \prod_{v \in E_{B}(\theta)} F_{v} \Big|
\exp \left( |\hbox{Im}\,p| \max_{\ell\in
L_{2}(\theta)} |x_{\ell}| \right)
\le F^{|E_{B}(\theta)|}
\prod_{v \in E_{B}(\theta) } {\rm e}^{-\x|\n_{v}|/2} ,
\label{eq:2.12} \end{equation}
provided $|\om|\sigma<\x/2$, and summability over the Fourier labels
in (\ref{eq:2.3}) is assured. The sum over $k$ in (\ref{eq:2.2})
produces a quantity bounded proportionally to the exponential
${\rm e}^{\Gamma|p|}$, for some positive constant $\Gamma$.
A comparison with \cite{GBD} shows that $\Gamma=1/\e_{0}$,
where $\e_{0}$ is the same as in the statement of the theorem.
In particular the Borel transform $x_{B}(t;p)$ of the series
(\ref{eq:2.2}) turns out to have an analytic continuation to the
strip $\Sigma_{\sigma}$, and admits there the bound
$|x_{B}(t;p)|\le C {\rm e}^{\Gamma|p|}$, for a suitable
constant $C$. Hence the integral
\begin{equation}
g(t;\e) := \int_{0}^{\infty} {\rm e}^{-p/\e} x_{B}(t;p) \, {\rm d} p 
\label{eq:2.13} \end{equation}
absolutely converges provided $0<\e<\e_{0}$.
So also the last two conditions for the formal series
of $x(t;\e)$ to be Borel summable are satisfied. 

That the solution $x(t;\e)$ describes a local attractor,
under the further condition $g'(c_{0})>0$,
follows from the analysis performed in \cite{BDG}.$\EP$

Note that, because of the analyticity properties of $x_{B}(t;p)$,
it follows, as a consequence of Nevanlinna's theorem \cite{S}, that
the function defined by the integral (\ref{eq:2.13}) is analytic in
the circle $C_{R}=\{ \e \in \CCC : {\rm Re}\,\e^{-1} > R^{-1}\}$,
with $R=\e_{0}$, and satifies the bound
\begin{equation} g(t;\e) =
\sum_{k=0}^{N-1} \e^{k} x^{(k)}(t) + \RRRR_{N}(\e) ,
\qquad |\RRRR_{N}(\e)| \le A B^{N} N! |\e|^{N} ,
\label{eq:2.14} \end{equation}
with constants $A$ and $B$ independent of $N$. This
is consistent with Proposition 5.3 of \cite{GBD}.

\zerarcounters
\section{Borel summability for $\boldsymbol{d}
\boldsymbol{=}\boldsymbol{2}$ and $\boldsymbol{\tau}
\boldsymbol{=}\boldsymbol{1}$}
\label{sec:3}

In the case of quasi-periodic forcing terms for $d=2$
we obtain the following result.

\begin{theorem}\label{thm:2}
Consider the system (\ref{eq:1.1}) for $d=2$, and assume that
$\oo$ satisfies the Diophantine condition (\ref{eq:1.3}) with $\tau=1$
and that the non-degeneracy condition (\ref{eq:1.4}) is fulfilled.
There exists $\e_{0}>0$ such that for $|\e|<\e_{0}$ there
is a quasi-periodic solution $x(t;\e)$ which has the same frequency
vector as the forcing term and is Borel summable at the origin.
If $g'(c_{0})>0$ such a solution describes a local attractor.
\end{theorem}

\prova Again we discuss explicitly the case $g(x)=x^{2}$
in (\ref{eq:1.1}). Let $\psi$ be a non-decreasing $C^{\infty}$
function defined in $\RRR_{+}$, such that
\begin{equation}
\psi(u) = \left\{
\begin{array}{ll}
1 \, , & \text{for } u \geq 1 \, , \\
0 \, , & \text{for } u \leq 1/2 \, ,
\end{array} \right.
\label{eq:3.1} \end{equation}
and set $\chi(u) := 1-\psi(u)$. Define, for all $n \in \ZZZ_{+}$,
$\chi_{n}(u) := \chi(2^{n}C_{0}^{-1}u/4)$ and $\psi_{n}(u) :=
\psi(2^{n}C_{0}^{-1}u/4)$.

With each line $\ell$ with zero momentum we associate a scale label
$n_{\ell}=-1$, while with each line witn non-zero momentum we
associate (arbitrarily) a scale label $n_{\ell} \in \ZZZ_{+}=
\{0\}\cup\NNN$. Then we can define cluster and
self-energy clusters as in \cite{G1,GBD}.
A cluster $T$ on scale $n$ is a maximal set of points and lines
connecting them such that all the lines have scales $n'\le n$
and there is at least one line with scale $n$. The lines entering
the cluster $T$ and the possible line coming out from it (unique if
existing at all) are called the external lines of the cluster $T$.
Given a cluster $T$ on scale $n$, we shall denote by $n_{T}=n$ the
scale of the cluster; we call $V(T)$, $E(T)$, $E_{W}(T)$, $E_{B}(T)$,
and $L(T)$ the set of vertices, of endpoints, of white endpoints,
of black endpoints, and of lines of $T$, respectively.
We call self-energy cluster any cluster $T$ such that
$T$ has only one entering line $\ell_{T}^{2}$ and one exiting
line $\ell_{T}^{1}$, and one has $\sum_{v\in E_{B}(T)} \nn_{v} = \vzero$.
With each line $\ell$ with momentum $\nn_{\ell}$ and scale $n_{\ell}$
we associate a renormalised propagator $g_{\ell}=
g^{[n_{\ell}]}(\oo\cdot\nn_{\ell};\e)$, still to be defined.
On the contrary the node factors are defined as in the previous
case (with the only trivial difference that now $\nn_{v}$,
replacing $\n_{v}$, is a $d$-dimensional vector).

Define the self-energy value $\VV_{T}(\oo\cdot\nn;\e)$ in terms
of the renormalised propagators and node factors as
\begin{equation}
\VV_{T}(\oo\cdot\nn;\e) = \Big( \prod_{\ell \in L(T)}
g^{[n_{\ell}]}(\oo\cdot\nn_{\ell};\e) \Big)
\Big( \prod_{v \in E(T) \cup V(T)} F_{v} \Big) ,
\label{eq:3.2} \end{equation}
where $\nn$ is the momentum of both the external lines of $T$.

We proceed as in Section VI of \cite{GBD}, with the only two differences
that we perform a preliminary summation by including
the contribution $-2 \e c_{0}$ (arising from the self-energy graphs
on scale $-1$) into the propagator $g^{[0]}(x;\e)$, and -- as
in the periodic case of Section \ref{sec:2} -- we associate
the factors $\e$ to the propagators with non-zero momentum.
Therefore we define\footnote{See footnote 1 in Section \ref{sec:2}.}
\begin{equation}
g^{[0]}(x;\e) = \frac{\e \psi_{0}(|x|) }{ ix (1+i\e x)-2\e c_{0}} ,
\qquad M^{[0]}(x;\e) = \e \sum_{k=1}^{\io} \sum_{T \in \SSS^{\RR}_{k,0}}
\VV_{T}(x;\e) ,
\label{eq:3.3} \end{equation}
whereas the propagators on scale $n\ge 1$ are defined as in \cite{GBD},
again with a factor $\e$ appearing in the numerator of the propagators
with non-zero momentum; this means that one has
\begin{equation} \begin{split}
g^{[n]}(x;\e) & = \frac{\e\chi_{0}(|x|) \ldots
\chi_{n-1}(|x|)\psi_{n}(|x|)}{ ix (1+i\e x) - \MM^{[n-1]}(x;\e)} , \\
\MM^{[n]}(x;\e) & = 
\sum_{p=1}^{n} \chi_{0}(|x|) \ldots \chi_{p-1}(|x|)\chi_{n}(|x|)
M^{[p]}(x;\e) , \qquad
M^{[n]}(x;\e) = \e \sum_{k=1}^{\io} \sum_{T \in \SSS^{\RR}_{k,n}}
\VV_{T}(x;\e) , 
\label{eq:3.4} \end{split} \end{equation}
where the set of renormalized self-energy clusters $\SSS^{\RR}_{k,n}$
is defined ad the set of self-energy clusters $T$ on scale $n_{T}=n$
and of order $k$ (that is with $|V(T)|+|E_{B}(T)|=k$). With respect
to \cite{GBD,G1} a further factor $\e$ appears in $M^{[n]}(x;\e)$,
$n \ge 0$, simply because there is one of such factors per node
(vertex or endpoint) with exiting line carrying a non-zero
momentum -- cf. Section 6 in \cite{G1} --, and we are associating
the factors $\e$ with the lines instead of the nodes.

An easy computation gives, for the Borel transform of $g^{[0]}(x;\e)$,
\begin{equation}
g^{[0]}_{B}(x;p) = \frac{\psi_{0}(|x|)}{ix} \exp
\left( -ip \left( x - 2 \frac{c_{0}}{x} \right) \right)
\qquad \Longrightarrow
\qquad \left| g^{[0]}_{B}(x;p) \right| \le
\frac{1}{|x|} {\rm e}^{(|x|+2|c_{0}|/|x|)|{\rm Im} \,p|} .
\label{eq:3.5} \end{equation}
If we set, for $n \ge 0$,
\begin{equation}
\widetilde g^{[n]}(x;\e) =
\frac{\e}{ ix (1+i\e x) - \MM^{[n-1]}(x;\e)}  \qquad
\forall |x| \le 2^{-(n-1)}C_{0} ,
\label{eq:3.6} \end{equation}
and define $M^{[n]}(x;\e)=\MM^{[n]}(x;\e)-
\MM^{[n-1]}(x;\e)$, we obtain the recursive equations
\begin{equation}
\left( \widetilde g^{[n]}(x;\e) \right)^{-1} =
\left( \widetilde g^{[n-1]}(x;\e) \right)^{-1} -
\chi_{0}(|x|) \ldots \chi_{n-1}(|x|)
\e^{-1} M^{[n-1]}(x;\e) \qquad n \ge 1 .
\label{eq:3.7} \end{equation}

By using these equations we can prove inductively the bound
\begin{equation}
\left| \widetilde g^{[n]}_{B}(x;p) \right| \le
\frac{K_{0}}{|x|} {\rm e}^{(c_{n}+c_{n}'|x|^{-1/2})|p|+
\ka_{0}|{\rm Im}\,p|(d_{n}|x|+d_{n}'|x|^{-1})} ,
\label{eq:3.8} \end{equation}
where $K_{0}$ and $\ka_{0}$ are two constants, and
the sequences $\{c_{n}\}_{n=0}^{\infty}$,
$\{c_{n}'\}_{n=0}^{\infty}$, $\{d_{n}\}_{n=0}^{\infty}$,
$\{d_{n}'\}_{n=0}^{\infty}$ are to be find out.

The proof proceeds as in Appendix A1 of \cite{CGGG}.
Set $x_{\ell}=\oo\cdot\nn_{\ell}$, and call $L_{0}(T)$ and
$L_{2}(T)$ the set of lines in $L(T)$ with zero momentum
and the set $L_{2}(T)=L(T)\setminus L_{0}(T)$, respectively.
First we use the inductive bound to obtain
\begin{equation} \begin{split}
\left| \left( \frac{M^{[N]}(x;\e)}{\e} \right)_{B} \right|
& \le \sum_{k=2}^{\infty} \sum_{T \in \SSS^{\RR}_{k,N-1}}
\Big( \prod_{\ell \in L_{0}(T)} |g_{\ell}| \Big)
\Big( \prod_{v \in E(T) \cup V(T)} |F_{v}| \Big) \\
& \qquad \Big( {\mathop{\prod}_{\ell \in L_{2}(T)}}^{\!\!\!\!*}
\frac{K_{0}}{|x_{\ell}|} {\rm e}^{(c_{n_{\ell}}+
c_{n_{\ell}}'|x_{\ell}|^{-1/2})|p|+ \ka_{0} (d_{n_{\ell}}|x_{\ell}|+
d_{n_{\ell}}'|x_{\ell}|^{-1})|{\rm Im}\,p|} \Big) \\
& \le \Big( \prod_{v \in E_{B}(\theta)} {\rm e}^{-\x|\n_{v}|} \Big)
\sum_{k=2}^{\infty} \Gamma^{k} \frac{|p|^{k-2}}{(k-2)!}
{\rm e}^{(c_{N-1}+c_{N-1}'2^{N/2})|p|+
\ka_{0} d_{N-1}'2^{N} |{\rm Im}\,p| } ,
\label{eq:3.9} \end{split} \end{equation} 
where $D_{0}=\Gamma^{2}$, $r_{N}=\Gamma + c_{N-1} + \Gamma_{0}
c_{N-1}'2^{N/2}$, for some $N$-independent constant $\Gamma_{0}$.
The bound in the last line of (\ref{eq:3.9}) has been obtained
by using part of the exponential decay (say one fourth)
of the node factors associated with the endpoints to control the
exponent $\ka_{0} d_{N-1}\max_{\ell\in L_{2}(T)}|x_{\ell}|$,
provided $d_{N-1}<d$ for some $N$-independent constant $d$ and
$|{\rm Im}\,p| \le \sigma$, with $\sigma$ small enough, more
precisely $\ka_{0} \sigma d |\om| < \xi/4$.

By explicitly performing the sum over $k$ we obtain from (\ref{eq:3.6})
\begin{equation}
\left| \left( \frac{M^{[N]}(x;\e)}{\e} \right)_{B} \right|
\le D_{0} {\rm e}^{r_{N}|p|} {\rm e}^{-\xi_{0} 2^{N}} ,
\label{eq:3.10} \end{equation} 
where we have used the bound $\sum_{v \in E_{B}(T)}|\nn_{v}| \ge
\Gamma_{1}2^{N}$, for a suitable constant $\Gamma_{1}$ -- see
formula (7.12) of \cite{G1} -- and again part of the exponential decay
(say another one fourth) of the node factors associated
with the endpoints to control the exponent $\ka_{0} d_{N-1}'2^{N}
|{\rm Im}\,p|$, provided again $d_{N-1}'<d'$ for some $N$-independent
constant $d'$ and $\ka_{0}d'\sigma < \xi \Gamma_{1}/4$;
in particular one finds $\xi_{0}=\Gamma_{1}\xi/4$.

Then, by using (\ref{eq:3.10}) and, once more,
the inductive bound, we obtain from (\ref{eq:3.7})
\begin{equation} \begin{split}
& \left| \widetilde g^{[N]}_{B}(x,p) \right| \le
\frac{K_{0}}{|x|} {\rm e}^{(c_{N-1}+c_{N-1}'|x|^{-1/2})|p|+
\ka_{0}(d_{N-1}|x|+d_{N-1}'|x|^{-1})|{\rm Im}\,p|} \\
& \qquad * \sum_{k=0}^{\infty}
\left( \left( D_{0} {\rm e}^{r_{N}|p|} {\rm e}^{-\xi_{0} 2^{N}}
\right) * \left( 
\frac{K_{0}}{|x|} {\rm e}^{(c_{N-1}+c_{N-1}'|x|^{-1/2})|p|+
\ka_{0}(d_{N-1}|x|+d_{N-1}'|x|^{-1})|{\rm Im}\,p|} \right)
\right)^{*k} ,
\label{eq:3.11} \end{split} \end{equation}
with $a^{*k}=a*a*\ldots*a$ ($k$ times). This gives
\begin{equation}
\left| \widetilde g^{[N]}_{B}(x,p) \right| \le
\frac{K_{0}}{|x|} \sum_{k=0}^{\infty} \frac{1}{(2k)!}\left(
\frac{K_{0}|p|^{2}}{|x|} D_{0} {\rm e}^{- \xi_{0}2^{N}} \right)^{k}
{\rm e}^{(r_{N}+c_{N-1}'|x|^{-1/2})|p|+
\ka_{0}(d_{N}|x|+d_{N}'|x|^{-1})|{\rm Im}\,p|} ,
\label{eq:3.12} \end{equation}
which implies the bound (\ref{eq:3.5}) for $n=N$, with $c_{N}=r_{N}=
\Gamma + c_{N-1} + \Gamma_{0} c_{N-1}' 2^{N/2}$, $c_{N}'=c_{N-1}'+
\sqrt{ K_{0} D_{0} {\rm e}^{- \xi_{0} 2^{N}}}$, $d_{N}=d_{N-1}$ and
$d_{N}'=d_{N-1}'$. In particular one has $d_{N}=d=1$ and $d_{N'}=d'=
2|c_{0}|$, so that there exists a constant $c>0$ such that $\max\{
c_{n}2^{-n/2},c_{n}',d_{n},d_{n'}\} \le c$ for all $n\ge 0$.

The bounds (\ref{eq:3.8}) for the Borel transforms of the
propagators can be used to obtain a bound on the
Borel transform $x_{B}(t;p)$ of $x(t;\e)$. We omit the details,
which can be derived exactly as in Appendix A1 of \cite{CGGG}.
Eventually one finds the bound
\begin{equation}
\left| x_{B}(t;p) \right| \le C_{1} {\rm e}^{C_{2}|p|^{2}} ,
\label{eq:3.13} \end{equation}
for suitable constants $C_{1}$ and $C_{2}$.
Again, the bound (\ref{eq:3.13}) and the analyticity properties
of $x_{B}(t;p)$ implies that $x(t;\e)$ is Borel summable,
and it can be written for $\e>0$ as
\begin{equation}
x(t;\e) = \int_{0}^{\infty} {\rm e}^{-p/\e} x_{B}(t;p)
\, {\rm d}p ,
\label{eq:3.14} \end{equation}
in terms of its Borel transform.

As in the case $d=1$ the last statement of the theorem
has been proved in \cite{BDG}.$\EP$

In the general case $g(x)\neq x^{2}$ in (\ref{eq:1.1})
the quantity $2c_{0}$ has to be replaced with $g'(c_{0})$, with
$g'(c_{0}) \neq 0$ by hypothesis. Then the discussion proceeds
as in Section VII of \cite{GBD}.

Note also that in the case $d=2$ and $\tau=1$
the Borel transform is still defined in a strip around the
real axis, but it does not satisfy any more an exponential 
bound like in the case $d=1$ (at least the argument given above
does not provide an estimate of this kind). Thus, we cannot apply
Nevanlinna's theorem to prove Borel summability.

\zerarcounters
\section{Bryuno frequency vectors}
\label{sec:4}

Let $\oo\in\RRR^{d}$ be a Bryuno vector. This means that
$B(\oo)<\infty$, with $B(\oo)$ defined in (\ref{eq:1.6}).

\begin{theorem} \label{thm:3}
Consider the system (\ref{eq:1.1}) for any $d\ge 2$, and assume that
$\oo$ satisfies the Bryuno condition $B(\oo)<\infty$ and
that the non-degeneracy condition (\ref{eq:1.4}) is fulfilled.
There exists $\e_{0}>0$ such that for all real $|\e|<\e_{0}$
there is a quasi-periodic solution with frequency vector $\oo$.
If $g'(c_{0})>0$ such a solution describes a local attractor.
\end{theorem}

For simplicity's sake we discuss the case $g(x)=x^{2}$
and $\e\in\RRR$, but the analysis can be easily generalised
to any analytic function $g$ (provided the non-degeneracy condition
(\ref{eq:1.4}) is satisfied). Furthermore the solution
can be showed to extend to a function analytic in $\e$
in the domain $\mathcal{C}_{R}$ defined in Section VI of \cite{GBD}
(cf. Figure 16 in \cite{GBD}).

Let $\psi(x)$ be the non-decreasing $C^{\infty}$ function
defined in (\ref{eq:3.1}) and set $\chi(x) := 1-\psi(x)$.
Define, for all $n \in \ZZZ_{+}$, $\chi_{n}(x) := \chi(\al_{n}^{-1}
(\oo) x/4)$ and $\psi_{n}(x) := \psi(\al_{n}^{-1}(\oo)x/4)$.

Set $g^{[-1]}(x;\e)=1$ and $M^{[-1]}(x;\e)=0$, and define iteratively
$g^{[n]}(x;\e)$ and $M^{[n]}(x;\e)$ as done in the case
of Diophantine vectors. This means that for $n=0$ we can define
$g^{[0]}(x;\e)$ and $M^{[0]}(x;\e)$ as in (\ref{eq:3.3}),
while for $n \ge 1$ we define
\begin{equation} \begin{split}
g^{[n]}(x;\e) & = \frac{\e\chi_{0}(|x|) \ldots \chi_{n-1}(|x|)
\psi_{n}(|x|)}{ix(1+i\e x)-\MM^{[n-1]}(x;\e)} , \\
\MM^{[n]}(x;\e) & = \sum_{p=0}^{n} \chi_{0}(|x|) \ldots
\chi_{p}(|x|) M^{[p]}(x;\e) , \qquad
M^{[n]}(x;\e) = \e \sum_{k=1}^{\infty}
\sum_{T \in \SSS^{\RR}_{k,n}} \VV_{T}(x;\e) ,
\label{eq:4.1} \end{split} \end{equation}
where $\SSS^{\RR}_{k,n}$ is the set of renormalised self-energy clusters
$T$ on scale $n$ and of order $k$, and the self-energy value
$\VV_{T}(x;\e)$ is defined as in (\ref{eq:3.2}).
Note that we are using the same definitions of Section \ref{sec:3},
in particular we are associating the factors $\e$ with the propagators
rather than with the nodes (contrary to what done in \cite{GBD}).
So far the only difference with respect to the case of the
standard Diophantine condition concerns the multiscale decomposition:
the factors $2^{n}C_{0}^{-1}$ appearing in $\chi_{n}$
and $\psi_{n}$ are substituted with $\al_{n}^{-1}(\oo)$.

\begin{lemma} \label{lem:1}
Assume that the renormalised propagators up
to scale $n-1$ can be bounded as
\begin{equation}
\left| g^{[n_{\ell}]}(\oo\cdot\nn_{\ell};\e) \right| \le C^{-1}
\al_{n_{\ell}}^{-\beta}(\oo)
\label{eq:4.2} \end{equation}
for some positive constants $\beta$ and $C$.
Then for all $p \le n-1$ the number $N_{p}(\theta)$ of lines
on scale $p$ in any renormalised tree $\theta$ and the number $N_{p}(T)$
of lines on scale $p$ in any renormalised self-energy cluster $T$  
are bounded both by
\begin{equation}
N_{p}(\theta) \le K 2^{-p} \sum_{v\in E_{B}(\theta)} |\nn_{v}| ,
\qquad N_{p}(T) \le K 2^{-p} \sum_{v\in E_{B}(T)} |\nn_{v}| ,
\label{eq:4.3} \end{equation}
for some positive constant $K$. If $|\e| < \e_{0}$,
with $\e_{0}$ small enough, then for all $p \le n-1$ one has
\begin{equation}
| M^{[p]}(x;\e) | \le D_{1} |\e|^{2} {\rm e}^{- D_{2} 2^{p}} , \qquad
| \partial_{x} M^{[p]}(x;\e) | \le
D_{1} |\e|^{2} {\rm e}^{- D_{2} 2^{p}} ,
\label{eq:4.4} \end{equation}
for some positive constants $D_{1}$ and $D_{2}$.
Only the constant $D_{1}$ depends on $\beta$. The constant $\e_{0}$
can be written as $\e_{0}= C_{1} \al_{n_{0}}^{\beta}$,
with $n_{0}=n_{0}(\oo,\beta)$ such that
\begin{equation}
K \beta \sum_{n=n_{0}+1}^{\infty} \frac{1}{2^{n}}
\log \frac{1}{\al_{n}(\oo)} \le \frac{\xi}{4} ,
\label{eq:4.5} \end{equation}
and $C_{1}$ a positive constant dependending on $C$ but not on $\beta$.
\end{lemma}

\prova The lemma can be proved by reasoning as in \cite{G1,G2}.
We simply sketch the proof, and omit the details.
First of all note that, if we define $n(\nn)=\{n\in\ZZZ_{+}\,:\, 2^{n-1}<
|\nn| \le 2^{n}\}$ then one has $|\oo\cdot\nn| \ge \al_{n(\nn)}(\oo)$.
Moreover $n'>n$ implies $\al_{n'}(\oo) \le \al_{n}(\oo)$, and
$\al_{n'}(\oo) < \al_{n}(\oo)$ implies $n'>n$.
Set $M(\theta) = \sum_{v\in E_{B}(\theta)} |\nn_{v}|$ and
$M(T) = \sum_{v\in E_{B}(T)} |\nn_{v}|$.
The bound on $N_{p}(\theta)$ is obtained by proving by induction
on the order of the renormalised tree that if $N_{p}(\theta) \neq 0$
then $N_{p}(\theta) \le 2 \, 2^{-p} M(\theta)-1$
Then, given a renormalised self-energy cluster $T\in \SSS^{\RR}_{k,n}$,
one proves first that $M(T) > 2^{n-1}$, hence, again by induction,
that if $N_{p}(T) \neq 0$ then $N_{p}(T) \le 2 \, 2^{-p} M(T)-1$.
Therefore (\ref{eq:4.3}) is proved. An important property is that
if a cluster $T$ has two external lines, with momenta $\nn$
and $\nn'$, respectively, with $\nn\neq\nn'$, both on scales greater
or equal to $n$, so that $|\oo\cdot\nn|\le \al_{n-1}(\oo)/4$
and $|\oo\cdot\nn'|\le \al_{n-1}(\oo)/4$, then one has
$|\oo\cdot(\nn-\nn')| < \al_{n-1}(\oo)$, hence $n(\nn-\nn') \ge n$,
so that $M(T) \ge |\nn-\nn'|>2^{n-1}$. For details we refer to \cite{G2}.

The bounds (\ref{eq:4.4}) are obtained by exploiting the just mentioned
bound on $M(T)$ and half the exponential decay factors
${\rm e}^{-\xi|\nn_{v}|}$ associated with the vertices and
endpoints internal to $T$ to derive the factors
${\rm e}^{- D_{2} 2^{p}}$, with $D_{2}$ independent of $\beta$,
and by using the fact that any self-energy cluster $T$
contributing to $M^{[p]}(x;\e)$ must be of order at least 2
to derive the factors $|\e|^{2}$.

Then for any $n_{0} \in \NNN$ and for any tree $\theta$,
we can bound each propagator on scale up to $n_{0}$
with $C^{-1} \al^{-\beta}_{n_{0}}(\oo)$ and the product of propagators
on scale greater than $n_{0}$ with
\begin{equation}
\prod_{n=n_{0}+1} \left( C^{-1}
\al^{-\beta}_{n}(\oo) \right)^{N_{n}(\theta)}
= C^{-\sum_{n=n_{0}+1}^{\infty} N_{n}(\theta)}
\exp \left( \beta M(\theta) \sum_{n=n_{0}+1}^{\infty}
\frac{1}{2^{n}} \log \frac{1}{\al_{n}(\oo)} \right) ,
\label{eq:4.6} \end{equation}
so that, by choosing $n_{0}$ according to (\ref{eq:4.5}),
the last exponential in (\ref{eq:4.6}) is controlled by
half the exponential decay factor ${\rm e}^{-\xi M(T)}$
arising from the node factors. Then the sum of the values
of all trees of order $k$ is bounded by
$(C^{-1} C' \al_{n_{0}}^{-\beta})^{k}$, for a suitable constant $C'$
-- taking into account all the constants other than $C$
and the sums over the trees. Hence also the assertion about
the dependence of $\e_{0}$ on $\al_{n_{0}}(\oo)$ follows,
and the proof of the lemma is complete.$\EP$

As in \cite{GBD} to prove existence of the quasi-periodic solution
we need the following result, which together with Lemma \ref{lem:1}
provides the proof of Theorem \ref{thm:3}.

\begin{lemma} \label{lem:2}
For real $\e$ small enough the renormalised propagators
satisfy the bounds (\ref{eq:4.2}) with $\beta=1$.
For $\e$ in the domain $\mathcal{C}_{R}$ in Figure 16 of \cite{GBD}
they satisfy the bounds (\ref{eq:4.2}) with $\beta=2$.
\end{lemma}

\prova The proof can be carried out exactly
as in \cite{GBD}. Indeed it is enough to show that
the propagators $g^{[n]}(x;\e)$ can be bounded proportionally
to $|x|^{-\beta}$, for $\e$ small enough in a suitable
domain, and this follows from Lemmata 6.2 to 6.5 of \cite{GBD},
independently on the particular Diophantine
condition assumed on $\oo$.$\EP$

The proof of the theorem is completed if we show that the
quasi-periodic solution is a local attractor if $g'(c_{0})>0$.
But this can be proved as in the case of Diophantine frequency vectors,
by reasoning as in \cite{BDG}: indeed the only property that we
need for the argument given in \cite{BDG} to work is the existence 
of the quasi-periodic solution.


\end{document}